\documentclass[10pt]{article}
\usepackage{amsmath}
\usepackage{amsfonts}
\usepackage{amssymb}

\begin{document}

\numberwithin{equation}{section}

\def\ye{y_{\varepsilon}}

\begin{center}
{\bf BOUNDARY VALUE PROBLEM FOR THE TIME-FRACTIONAL TELEGRAPH EQUATION WITH CAPUTO DERIVATIVES}
\end{center}
\begin{center}
{\bf M.O. Mamchuev }
\end{center}

In this paper the Green formula for the operator of fractional differentiation in Caputo sense is proved.
By using this formula the integral representation of all regular in a rectangular domains solutions is obtained 
in the form of the Green formula for operator generating the time-fractional telegraph equation.
The unique solutions of the initial-boundary value problem with boundary conditions of first kind is constructed.

The proposed approach can be used to study the more general evolution FPDE as well as ODE with Caputo derivatives.

 \medskip


 \smallskip


\section{Introduction}\label{sec:1}

\setcounter{section}{1}
\setcounter{equation}{0}

Consider the equation
\begin{equation}\label{mam1}
{\bf L} u(x,y)\equiv \partial_{0y}^{\alpha}u(x,y)+b\partial_{0y}^{\beta}u(x,y)
-\frac{\partial^2}{\partial x^2}u(x,y)+cu(x,y)=f(x,y),
\end{equation}
where  $\alpha=2\beta\in]0,2[,$
$b$ and $c$ are given real numbers, 
$f(x,y)$ is a given real function,
$\partial_{0y}^{\nu}$ is the Caputo fractional differentiation operator of order $\nu.$ 

Fractional partial differential equations 
find their application in various fields of science, such as physics, 
chemistry, biology, economics, sociology, etc. See for example the monographs \cite{n}, \cite{Hilfer}, \cite{Uchaikin}, 
and references therein.
Equations of (\ref{mam1}) type are used to describe anomalous diffusion processes observed 
in experiments related to blood circulation \cite{Cascaval-2002},
iterated Brownian motion, and telegraph processes with Brownian time
\cite{Orsinger-2003}, \cite{Orsinger-2004}, model of diffusion-drift charge carrier transport in layers with a fractal structure \cite{Rekh2016}.
Note that the fractional advection-dispersion equation with constant coefficient,  
which has been studied in papers  \cite{Liu-2003},  \cite{Huang-2005},
can be reduced to equation (\ref{mam1}), by changing of unknown function. 

Note also that the boundary value problems for the analogue of equation (\ref{mam1}) with Rimman-Liouville derivatives
with $0<\alpha=2\beta<1$ has been studied in \cite{izvkbnc04} and \cite{daman05}, and in 
\cite{du15-5} and \cite{du15-9} in general case, when $0<\alpha=2\beta<2$.

For an extensive  bibliography on this subject see paper \cite{Mamchu2016arx} and references therein.
 
For the operator of fractional differentiation in Rimman-Liouville sense Green formula was proved
in papers \cite {Ps-2003}.
Green's function method of solving the boundary value problems for the diffusion-wave equation 
was implemented 
on the basis of this formula 
\cite {Ps-2003}, \cite {PsMon}.

We obtain Green formula for Caputo fractional differentiation operator.
This formula allows us to implement Green function method for equation with Caputo derivatives. 
This is shown by the time-fractional telegraph equation (\ref{mam1}).
The general representations of solutions is obtained and the theorems of solution existence and uniquiness for boundary value problem is proved.

\section{Preliminaries}\label{sec:2}

\setcounter{section}{2}
\setcounter{equation}{0}

The operator of fractional integro-differentiation in Riemman-Liouville sense $D_{ay}^{\nu}$ 
of order $\nu$ is defined by the formula \cite [p. 9]{n}: 
$$
D_{ay}^{\nu}g(y)=\left \{
\begin{array}{ll}
\frac{{\rm sign}(y-a)}{\Gamma(-\nu)}\int\limits_a^yg(s)|y-s|^{-\nu-1}ds, &
\nu< 0,\\
g(y), & \nu= 0,\\
{\rm sign}^n(y-a)\frac{d^n}{dy^n}D_{ay}^{\nu-n}g(s),\,  &
n-1<\nu\leq n, \, n\in \mathbb{N},
\end{array}
\right.$$
where $\Gamma(z)$ is a Eullers gamma-function.

The fractional Caputo derivative \cite{Caputo} 
is determined by relation   \cite [p. 11]{n}
$$\partial_{ay}^{\nu}g(y)={\rm{sign}}^n(y-a)D_{ay}^{\nu-n}g^{(n)}(y), \quad
n-1<\nu\leq n, \quad n\in \mathbb{N}.$$

We have the following analogue of the Newton-Leibniz formula in the fractional calculus
\cite[p. 11]{n}
\begin{equation} \label{vfnl}
D^{\delta}_{ay}D^{\nu}_{ay}g(y)=D^{\delta+\nu}_{ay}g(y)-
\sum\limits_{k=1}^{n}\frac{|y-a|^{-\delta-k}}{\Gamma(1-\delta-k)}
\lim\limits_{y\rightarrow a}D^{\nu-k}_{ay}g(y),
\end{equation}
where $n-1<\nu\leq n,$ $n\in \mathbb{N}.$
In case $\delta=-\alpha,$ $\nu=\alpha,$ $n-1<\alpha\leq n$ the formula (\ref{vfnl})
has the form
\begin{equation}\label{mam33}
D^{-\alpha}_{ay}D^{\alpha}_{ay}g(y)=g(y)-
\sum\limits_{k=1}^{n}\frac{|y-a|^{\alpha-k}}{\Gamma(1+\alpha-k)}
\lim\limits_{y\rightarrow a}D^{\alpha-k}_{ay}g(y),
\end{equation}

The formula 
\begin{equation} \label{vfh}
\int\limits_a^{y}f(s)D_{a s}^{\nu}h(s)ds=
\int\limits_{a}^{y}h(s)D_{ys}^{\nu}f(s)ds, \quad \nu<0,
\end{equation}
is known as a formula of fractional integration by parts \cite[p. 34]{n}.

The following formula of fractional intero-differentiation of power function
\begin{equation}\label{anm1}
D_{ay}^{\alpha}|y-a|^{\mu-1}=
\frac{\Gamma(\mu)}{\Gamma(\mu-\alpha)}|y-a|^{\mu-\alpha-1}.
\end{equation}
holds for each $\mu>0$ and $\alpha\in {\Bbb R}.$ 
 
Here is the definition of the operator $S_{\eta y}^{\delta},$  acting on a function  $g(s)$ 
by the formula \cite [p. 22]{n}
\begin{equation} \label{vsd}
S_{sy}^{\delta}g(\xi)=\frac{{\rm sign}(y-s)}{\pi}
\int\limits_{s}^y
\frac{g(\tau)}{\tau-\xi}
\left|\frac{\tau-s}{s-\xi}\right|^{\delta}d\tau.
\end{equation}

The following two assertions are hold \cite[p. 117, 118]{PsMon}.

{\bf Lemma 2.1.} 
{\it Let $0\leq \eta<y_{\varepsilon}<y,$ $y_{\varepsilon}=y-\varepsilon.$ 
Then the estimate 
$$\left|S_{y_{\varepsilon}y}^{\delta}g(\eta)\right|
\leq C\varepsilon^{\gamma+\delta-\theta}(y_{\varepsilon}-\eta)^{\theta-\delta-1},
\quad \theta \in[0,1]. $$
holds for any function $g(\eta)\leq C(y-\eta)^{\gamma-1},$  and $0<\gamma<1.$  
}

{\bf Lemma 2.2.} 
{\it Let $y>y_{\varepsilon}>\eta,$ $y_{\varepsilon}=y-\varepsilon,$ $n-1<\alpha\leq n,$ $n\in \mathbb{N}.$
Then the relation  
\begin{equation}\label{mam31}
D_{\ye\eta}^{-\alpha}D_{y\eta}^{\alpha}g(\eta)=
g(\eta)-\sum\limits_{k=1}^{n}\frac{(\ye-\eta)^{\alpha-k}}{\Gamma(\alpha-k+1)}
D_{y\ye}^{\alpha-k}g(y_{\varepsilon})-
\gamma_n^{\alpha}S_{\ye y}^{n-\alpha}g(\eta),
\end{equation}
where $\gamma_n^{\alpha}=\sin  \pi(\alpha-n),$ 
holds for any function $g(\eta)$ continuous in the neighborhood of the point $\eta.$
}

The hipergeometric function is determined by the relation
$$F(a,b,c;z)=\sum\limits_{k=0}^{\infty}\frac{(a)_k(b)_k}{(c)_kk!}z^k,
\quad a,b,c\in {\mathbb C}, \quad c\not=0,-1,..., \quad |z|<1,$$
where $(g)_n=\Gamma(g+n)/\Gamma(g).$

Following assertions hold
\begin{equation}\label{hgf1}
F(a,b,c;z)=\frac{1}{B(b,c-b)}\int\limits_0^1t^{b-1}(1-t)^{c-b-1}(1-zt)^{-a}dt,
\end{equation}
for ${\rm Re}\, c >{\rm Re}\, b>0,$ $|z|<1,$
\begin{equation}\label{hgf2}
F(a,b,a;z)=(1-z)^{-b}.
\end{equation}

Let us consider the function
$$\Gamma(x,y)=\frac{1}{2}\int\limits_{|x|}^{\infty}
h_0(x,\tau)g(y,\tau)d\tau,$$
where
$$
h_0(x,\tau)=_0F_1\left[ 1; a(\tau^2-x^2)\right],
\quad
g(y,\tau)=\frac{e^{b_1\tau}}{y}\phi(-\beta,0;-\tau y^{-\beta}),
$$                                                   
$_0F_1(\nu;z)=\sum_{k=0}^{\infty}\frac{\Gamma(\nu)}{\Gamma(k+\nu)}\frac{z^k}{k!}$  
is the confluent hypergeometric limit function \cite [p. 32]{Kilbas-2006},
$\phi(\delta,\mu;z)=\sum_{k=0}^{\infty}\frac{1}{\Gamma(\delta k+\mu)}\frac{z^k}{k!}$  is the Wright function  \cite {Wright-1933}, 
$b_1=-b/2,$ and  $a=b_1^2-c.$

The following properties of the function $\Gamma(x,y)$ are hold \cite{du15-5}, \cite{Mamchu2016arx}.

{\bf Lemma 2.3.}
{\it
The relations  
$$D^{\delta}_{0y}D^{\nu}_{0y}\Gamma(x,y)=D^{\nu}_{0y}D^{\delta}_{0y}\Gamma(x,y)
=D^{\delta+\nu}_{0y}\Gamma(x,y)$$
hold for arbitrary $\delta,\nu \in \mathbb{R}.$ 
}

{\bf Lemma 2.4.}
{\it
The relation holds 
$$2\frac{\partial^m}{\partial x^m}D_{0y}^{\nu}\Gamma(x,y)=
\int\limits_{|x|}^{\infty}D_{0y}^{\nu-\beta k}g(y,\tau)
\frac{\partial^m}{\partial x^m}L_1^kh_0(x,\tau)d\tau-
$$
$$- \mathop{\rm sgn}(m)\mathop{\rm sgn}(x) \sum\limits_{j=1}^{m}
\frac{\partial^{j-1}}{\partial x^{j-1}}\left[D_{0y}^{\nu-\beta k}g(y,|x|)
\left(\frac{\partial^{m-j}}{\partial x^{m-j}}L_1^kh_0(x,\tau)\right)\big|_{\tau=|x|}
\right]+                                             
$$
$$+\mathop{\rm sgn}(k)\sum\limits_{i=1}^{k}\frac{\partial^m}{\partial x^m}
\left[D_{0y}^{\nu-\beta i}g(y,|x|)
\left(L_1^{i-1}h_0(x,\tau)\right)\big|_{\tau=|x|}\right],
$$
where  $k=0$ for $\nu\leq 0, $ and $k$ such that $\beta (k-1)\nu\leq \beta k$ $(k\in \mathbb{N})$
for $\nu>0,$
is satisfied for arbitrary $m\in \mathbb{N}\cup \{0\}$ and $\nu\in \mathbb{R},$
and $L_1=\frac{\partial}{\partial \tau}+b_1$
is the differential operator of the first order.
}

{\bf Lemma 2.5.}
{\it
The estimate 
\begin{equation} \label{kdxdyg}
\left|\frac{\partial^m}{\partial x^m}D_{0y}^{\nu}\Gamma(x,y)\right| 
\leq C |x|^{-\theta}y^{\beta(1-m+\theta)-\nu-1}, 
\quad \theta\geq 0, 
\end{equation}
where $C$ is a positive constant, holds for arbitrary
 $m\in \mathbb{N}\cup \{0\}$ and $\nu\in \mathbb{R}.$ 
}

We denote
$$\Omega=\{(x,y):a_1<x<a_2, 0<y<T\}, \quad
\Omega_y=\{(\xi,\eta):\;a_1<\xi<a_2,\;0<\eta<y\}.$$

{\bf Lemma 2.6.}
{\it
The relation ${\bf L}\frac{\partial^{m}}{\partial x^{m}}\Gamma(x,y)=0$
holds for $m=0,1,$ and for all $(x,y)\in\Omega.$
}

{\bf Lemma 2.7.}
{\it
The relation 
$${\bf L^*}\Gamma(x-\xi,y-\eta)\equiv 
\left(D_{y\eta}^{\alpha} +bD_{y\eta}^{\beta}-\frac{\partial^2}{\partial \xi^2}+c\right)\Gamma(x-\xi,y-\eta)=0 $$
holds for all $(\xi,\eta)\in\Omega_y$ and for fixed $(x,y)\in \Omega.$ 
}

{\bf Lemma 2.8.}
{\it
The relation
$$\lim\limits_{s\rightarrow y}\int\limits_{x_1}^{x_2}q(\xi)D^{\alpha-1}_{y\eta}
\Gamma(x-\xi,y-\eta)d\xi=q(x), \quad x_1<x<x_2.
$$
holds for each function $g(x)\in C[x_1,x_2],$ $a_1\le x_1<x_2\le a_2,$ 
}

{\bf Lemma 2.9.}
{\it
The relation
$$\lim\limits_{\xi\to x\pm 0}\int\limits_{\delta}^yp(\eta)
\frac{\partial}{\partial x}\Gamma(x-\xi,y-\eta)d\eta=\mp \frac{1}{2}p(y), \quad 0\leq\delta< y,
$$
holds for each function $p(y)\in C[\delta,T].$ 
}

By $C^{1,q}[a_1, a_2]$ we denote the space of continuously
differentiable functions on $[a_1, a_2]$ the first derivatives of which  
satisfy the H\"older condition with exponent $q.$

{\bf Lemma 2.10.} 
{\it
Let the functions $\tau_k(x)$ $(k=1,n)$  satisfy the conditions
\begin{equation} \label{equstk}
\begin{array}{c}
\tau_n(x)\in C[a_1,a_2]; \\
\tau_1(x)\in C^{1,q}[a_1,a_2],\quad q>\frac{1-\beta}{\beta}, 
\quad \mbox{for }\, n=2.
\end{array}
\end{equation} 
Then function 
$$u_0(x,y)=\sum\limits_{k=1}^{n}\int\limits_{a_1}^{a_2}
\tau_k(\xi)\left[D_{0y}^{\alpha-k}+(2-k)bD_{0y}^{\beta-k}\right]\Gamma(x-\xi,y)dt 
$$
is a solution of the equation (\ref{mam1}) in the class 
$\partial_{0y}^{\alpha}u_0,$ $\partial_{0y}^{\beta}u_0,$ 
$\frac{\partial^2}{\partial x^2}u_0\in C(\Omega),$ 
and satisfies the conditions  
\begin{equation} \label{neq2220}
\lim\limits_{y \rightarrow 0} 
\frac{\partial^{k-1}}{\partial y^{k-1}} u_0(x,y)=\tau_k (x),
\quad a_1<x<a_2,  \quad k=1,n. 
\end{equation} 
}

{\bf Lemma 2.11.} 
{\it
Let 
$f(x,y)\in C(\bar\Omega),$
$f_y(x,y)\in C(\Omega)\cup L(\Omega).$
Then the function 
$$u_f(x,y)=\int\limits_0^y \int\limits_{a_1}^{a_2}\Gamma(x-\xi,y-\eta)f(\xi,\eta)d\xi d\eta$$ 
is a regular solution of equation} (\ref{mam1}) {\it in domain $\Omega$ satisfying 
the homogenous condition} (\ref{neq2220}).

\section{Problem statement and main results}\label{sec:3}

\setcounter{section}{3}
\setcounter{equation}{0}
Set $J=\{(x,y):a_1<x<a_2, y=0\};$
and $n\in \{1,2\}$ is a number such that $n-1<\alpha\leq n.$

A regular solution of equation (\ref{mam1}) in domain $\Omega$ is defined as a  
function $u=u(x,y)$ of the class 
$u\in C(\overline{\Omega}),$ 
$\frac{\partial^{n-1} u}{\partial y^{n-1}}\in C(\Omega\cup J),$ 
$\partial_{0y}^{\alpha}u,\, \partial_{0y}^{\beta}u, 
\,\frac{\partial^2}{\partial x^2}u \in C(\Omega),$
satisfying equation (\ref{mam1}) at all points $(x,y)\in \Omega.$

{\bf Theorem 3.1.}
{\it
Let $f(x,y)\in C(\overline{\Omega}),$ 
and $\tau_k(x)\in C[a_1,a_2],$ $k=1,n.$
Then any regular solution $u(x,y)$ of equation} (\ref{mam1}) {\it in the domain $\Omega$
satisfying the condition
\begin{equation} \label{oos121} 
\lim\limits_{y\to 0}\frac{\partial^{k-1}}{\partial y^{k-1}}u(x,y)=\tau_k(x), 
\quad a_1< x<a_2, \quad k=\overline{1,n}, 
\end{equation}
and such that $u_x\in C(a_1\leq x\leq a_2, 0<y<T),$ $u_x(a_1,y),$ $u_x(a_2,y)\in L[0,T],$
can be represented in the form
$$u(x,y)=\sum\limits_{k=1}^{n}\int\limits_{a_1}^{a_2}
\tau_k(\xi)\left[D_{0y}^{\alpha-k}+(2-k)bD_{0y}^{\beta-k}\right]G(x,y;\xi,0)d\xi+$$
$$+\sum\limits_{i=1}^{2}(-1)^{i}\int\limits_{0}^{y}
[G(x,y;a_i,\eta)u_{\xi}(a_i,\eta)-G_{\xi}(x,y;a_i,\eta)u(a_i,\eta)]d\eta+$$
\begin{equation} \label{kkkkk}
+\int\limits_{a_1}^{a_2}\int\limits_0^y
[G(x,y;\xi,\eta)f(\xi,\eta)-u(\xi,\eta)h(x,y;\xi,\eta)]d\xi d\eta,
\end{equation}
where
$\eta^{\alpha-n}h(x,y;\xi,\eta)\in L(\Omega_y),$
$G(x,y;\xi,\eta)=\Gamma(x-\xi,y-\eta)-V(x,y;\xi,\eta),$ 
and
$V\equiv V(x,y;\xi,\eta) \equiv V(x,\xi,y-\eta)$ 
is a solution of the equation 
\begin{equation} \label{lvh}
{\bf L}^*V(x,y;\xi,\eta)=h(x,y;\xi,\eta)
\end{equation}
in the class                                                            
$\frac{\partial^2}{\partial t^2}V\in C(\Omega\times\Omega_y)\cup L(\Omega_y),$ 
$V, \frac{\partial}{\partial \xi}V, D_{y\eta}^{\alpha}V, D_{y\eta}^{\beta}V
\in C(\Omega\times\overline{\Omega}_y),$
with the condition
\begin{equation} \label{oos12}
\lim\limits_{\eta\to y}D_{y\eta}^{\beta k-i}V(x,y;\xi,\eta)=0, \quad
k=1,n, \quad i=1,k.
\end{equation}
}

{\bf Boundary value problem.} 
In the domain $\Omega$ find solution 
$u(x,y)$ of equation (\ref{mam1})  with the conditions  (\ref{oos121})  and
\begin{equation} \label {zad00}
u(a_1,y)=\varphi_1(y), \quad  
u(a_2,y)=\varphi_2(y), \quad  0<y<T,
\end{equation}                 
where the $\tau_k(x),$ $(k=1,n),$ $\varphi_1(y)$ and $\varphi_2(y)$ are given functions.

The function $G(x,y;\xi,\eta)=\Gamma(x-\xi,y-\eta)-V(x,y;\xi,\eta),$ which satisfies the conditions
$$
\lim\limits_{t \rightarrow a_1}G(x,y;\xi,\eta)=0,\ 
\lim\limits_{t \rightarrow a_2}G(x,y;\xi,\eta)=0, \ y\not=\eta,
$$
where  $V\equiv V(x,y;\xi,\eta)$ is a function  under condition of Theorem 3.1
with $h(x,y;\xi,\eta)\equiv 0,$
we call {\it a Green function of the problem}  (\ref{mam1}), (\ref{oos121}), (\ref{zad00}).

{\bf Theorem 3.2.}
{\it
The function 
$$G(x,y;\xi,\eta)=\sum\limits_{m=-\infty}^{+\infty}
\bigg[\Gamma(X_1^m,y-\eta)-\Gamma(X_2^m,y-\eta)\bigg],$$
where 
$$X_1^m=2m(a_2-a_1)+x-\xi, \quad X_2^m=2m(a_2-a_1)+x+\xi-2a_1,$$
is the Green function of problem} 
(\ref{mam1}), (\ref{oos121}), (\ref{zad00}).

{\bf Theorem 3.3.}
{\it
Let the functions 
$\tau_k(x)$ $(k=1,n)$ satisfy the conditions} (\ref{equstk}),
{\it let $\varphi_i(y)\in C[0,T],$ $(i=1,2),$ and $f(x,y)\in C(\bar\Omega),$
\begin{equation} \label{mam4}
f_y(x,y)\in C(\Omega)\cup L(\Omega),
\end{equation}
and the relations 
\begin{equation} \label{usog}
\varphi_i(0)=\tau_1(a_i), \quad i=1,2
\end{equation}
hold.
Then the function
$$u(x,y)=\sum\limits_{k=1}^{n}\int\limits_{a_1}^{a_2}
\tau_k(\xi)\left[D_{0y}^{\alpha-k}+(2-k)bD_{0y}^{\beta-k}\right]G(x,y;\xi,0)d\xi+ $$
$$+\int\limits_{0}^{y}
\left[\frac{\partial}{\partial \xi}G(x,y;a_1,\eta)\varphi_1(\eta)-
\frac{\partial}{\partial \xi}G(x,y;a_2,)\varphi_2(\eta)\right]d\eta+$$
\begin{equation} \label{sp00}
+\int\limits_{0}^y\int\limits_{a_1}^{a_2} G(x,y;\xi,\eta)f(\xi,\eta)d\xi d\eta.
\end{equation}
is a unique regular solution of the problem} (\ref{mam1}), (\ref{oos121}), (\ref{zad00}).

{\bf Remark.}
If $ 2/3 <\alpha \leq 1$ it is sufficient 
that the function $ f(x, y) $ satisfies H\"older condition with respect to the 
variable $ x $ with exponent
$q>\frac{1-2\beta}{\beta},$ instead of the condition (\ref{mam4}) of Theorem 3.3.

\section{General representation of the solution }
\label{sec:4}

\setcounter{section}{4}
\setcounter{equation}{0}

\subsection{Green formula for the Caputo operator}
\label{subsec2}

\noindent

Following equality is true 
\begin{equation}\label{fnlc}
D_{ay}^{-\alpha+m-1}\partial_{ay}^{\alpha}g(y)={\rm sgn}^{m-1}(y-a) g^{(m-1)}(y)
-\sum\limits_{i=m}^{n}\frac{(y-a)^{i-m}}{(i-m)!}g^{(i-1)}(a).
\end{equation}

Indeed, using formula (\ref{vfnl}) we have
$$D_{ay}^{\delta}\partial_{ay}^{\nu}g(y)=D_{ay}^{\delta+\nu}g(y)-\sum\limits_{k=1}^{n}\frac{(y-a)^{n-\delta-\nu-k}}{\Gamma(n-\delta-\nu-k+1)}g^{(n-k)}(a).$$
By setting in the last relation $\delta=-\alpha+m-1,$   $\nu=\alpha,$ 
and by using the relations $\frac{1}{\Gamma(-p)}=0,$ ${\Gamma(p+1)=p!},$ 
which hold for $p\in {\mathbb N},$ we obtain  
$$D_{ay}^{-\alpha+m-1}\partial_{ay}^{\alpha}g(y)={\rm sgn}^{m-1}(y-a) g^{(m-1)}(y)
-\sum\limits_{k=1}^{n-m+1}\frac{(y-a)^{n-m-k+1}}{(n-m-k+1)!}g^{(n-k)}(a).$$
Replacing $i=n-k+1$ in the last relation get (\ref{fnlc}).

\noindent
{\bf Lemma 4.1.} 
{\it Let $\eta<\ye<y<T,$ $y_{\varepsilon}=y-\varepsilon,$ $n-1<\alpha\leq n,$ $n\in {\mathbb N}.$
Then the relations 
$$\int\limits_{0}^{\ye}
\left(v\partial_{0\eta}^{\alpha}u-uD_{y\eta}^{\alpha}v\right)d\eta
\sum\limits_{k=1}^{n}u^{(k-1)}(\ye)D_{y\ye}^{\alpha-k}v(y;\ye)-$$
\begin{equation}\label{fgco}
-\sum\limits_{k=1}^{n}u^{(k-1)}(0)D_{y\eta}^{\alpha-k}v(y;\eta)\big |_{\eta=0}+
\gamma_n^{\alpha}
\int\limits_{0}^{\ye}\partial_{0\eta}^{\alpha}u \cdot S_{\ye y}^{n-\alpha}v d\eta.
\end{equation}
hold for the function $u(\eta)$ and $v(y;\eta)$  such that 
$u(\eta)\in C[0,T],$ $u^{(n)}(\eta)\in L[0,y],$
$v(y;\eta)\in L([0,T]\times [0,y])\cap C([0,T]\times[0,\ye]).$
}

\noindent
{\bf Proof.}
Use the relation (\ref{mam31}), (\ref{vfh}) and definition of fractional integral operator, we can transform the next integral
$$\int\limits_{0}^{\ye}v(y;\eta)\partial_{0\eta}^{\alpha}u(\eta)d\eta=$$
$$=\int\limits_{0}^{\ye}\partial_{0\eta}^{\alpha}u(\eta)
\left[D_{\ye\eta}^{-\alpha}D_{y\eta}^{\alpha}v(y;\eta)
+\sum\limits_{k=1}^{n}\frac{(\ye-\eta)^{\alpha-k}}{\Gamma(\alpha-k+1)}D_{y\ye}^{\alpha-k}v(y;\ye)\right]d\eta+
$$
$$+\gamma_n^{\alpha}   	
\int\limits_{0}^{\ye}\partial_{0\eta}^{\alpha}u(\eta) S_{\ye y}^{n-\alpha}v(y;\eta) d\eta=
$$
$$
=\int\limits_{0}^{\ye}D_{0\eta}^{-\alpha}\partial_{0\eta}^{\alpha}u(\eta) \cdot 
D_{y\eta}^{\alpha}v(y;\eta) d\eta
+\sum\limits_{k=1}^{n}D_{y\ye}^{\alpha-k}v(y;\ye)\cdot D_{0\ye}^{-\alpha+k-1}\partial_{0\ye}^{\alpha}u(\ye) +
$$
\begin{equation}\label{eql31}
+\gamma_n^{\alpha} 	
\int\limits_{0}^{\ye}\partial_{0\eta}^{\alpha}u(\eta) 
S_{\ye y}^{n-\alpha}v(y;\eta) d\eta.
\end{equation}

Apply the formula (\ref{fnlc}) with $m=1$ and $a=0$ and integration by parts, for the first summand 
in the right hand side of last relation we obtain
$$
\int\limits_{0}^{\ye}D_{0\eta}^{-\alpha}\partial_{0\eta}^{\alpha}u(\eta) \cdot 
D_{y\eta}^{\alpha}v(y;\eta) d\eta= $$
$$=\int\limits_{0}^{\ye}u(\eta)D_{y\eta}^{\alpha}v(y;\eta)d\eta-
\sum\limits_{k=1}^{n}\frac{u^{(k-1)}(0)}{(k-1)!}\int\limits_{0}^{\ye}\eta^{k-1}D_{y\eta}^{\alpha}v(y;\eta)d\eta=
$$
$$=\int\limits_{0}^{\ye}u(\eta)D_{y\eta}^{\alpha}v(y;\eta)d\eta-
\sum\limits_{k=1}^{n}u^{(k-1)}(0)D_{y \eta}^{\alpha-k}v(y;\eta)\big |_{\eta=0}+ $$
\begin{equation}\label{eql32}
+\sum\limits_{k=1}^{n}u^{(k-1)}(0)\sum\limits_{j=1}^{k}
\frac{\ye^{k-j}}{(k-j)!}D_{y\ye}^{\alpha-j}v(y;\ye).
\end{equation}

And with the help of the formula (\ref{fnlc}), for the second summand we obtain
$$
\sum\limits_{k=1}^{n}D_{y\ye}^{\alpha-k}v(y;\ye)D_{0\ye}^{-\alpha+k-1}\partial_{0\ye}^{\alpha}u(\eta) = $$
\begin{equation}\label{eql33}
=\sum\limits_{k=1}^{n}u^{(k-1)}(\ye)D_{y\ye}^{\alpha-k}v(y;\ye)-\sum\limits_{k=1}^{n}D_{y\ye}^{\alpha-k}v(y;\ye)\sum\limits_{j=k}^{n}u^{(j-1)}(0)\frac{\ye^{j-k}}{(j-k)!}.
\end{equation}

From the relations (\ref{eql31}), (\ref{eql32}), (\ref{eql33}) we obtain  (\ref{fgco}).
The proof of Lemma 4.1 is complete.

\noindent
{\bf Lemma 4.2.} 
{\it Let $y>y_{\varepsilon}>\eta,$ $y_{\varepsilon}=y-\varepsilon,$ $n-1<\alpha\leq n,$ $n\in {\mathbb N}.$
Then the relations  
\begin{equation}\label{mam331}
D_{\ye\eta}^{\alpha-k}S_{\ye y}^{n-\alpha}g(\eta)=
(-1)^{n-k}\frac{\Gamma(\alpha-k+1)}{\pi} \int\limits_{\ye}^{y}\frac{g(\tau)d\tau}{(\tau-\eta)^{\alpha-(k-1)}},  \quad k=0,1,...,n,
\end{equation}
hold for any function $g(\eta)\in L[\ye,y].$
}

\noindent
{\bf Proof.}
By changing the order of integration, in case 
when $k=n,$ we obtain
$$D_{\ye\eta}^{\alpha-n}S_{\ye y}^{n-\alpha}g(\eta)=
\frac{1}{\pi\Gamma(n-\alpha)}\int\limits_{\eta}^{\ye}(\xi-\eta)^{-\alpha+n-1}d\xi
\int\limits_{\ye}^{y}
\frac{g(\tau)}{\tau-\xi}\left(\frac{\tau-\ye}{\ye-\xi}\right)^{n-\alpha}d\tau=$$
\begin{equation}\label{mam333}
=\frac{1}{\pi\Gamma(n-\alpha)}
\int\limits_{\ye}^{y} g(\tau)(\tau-\ye)^{n-\alpha}d\tau
\int\limits_{\eta}^{\ye}(\xi-\eta)^{-\alpha+n-1}(\tau-\xi)^{-1}(\ye-\xi)^{\alpha-n}d\xi.
\end{equation}

By taking into account the formulas (\ref{hgf1}), and (\ref{hgf2}), and replacing the integration variables 
$\xi=(\ye-\eta)z+\eta,$ we can transform the following integral
$$I_1(\ye,\eta,\tau)=
\int\limits_{\eta}^{\ye}(\xi-\eta)^{-\alpha+n-1}(\tau-\xi)^{-1}(\ye-\xi)^{\alpha-n}d\xi=$$
$$=(\tau-\eta)^{-1}\int\limits_{0}^{1}z^{-\alpha+n-1}(1-z)^{\alpha-n}
\left(1-\frac{\ye-\eta}{\tau-\eta}z\right)^{-1}dz=$$
$$=\Gamma(n-\alpha)\Gamma(\alpha-n+1)(\tau-\eta)^{-1}
F\left(1,n-\alpha,1;\frac{\ye-\eta}{\tau-\eta}\right)=
$$
\begin{equation}\label{mam334}
=\Gamma(n-\alpha)\Gamma(\alpha-n+1)(\tau-\eta)^{n-\alpha-1}(\tau-\ye)^{\alpha-n}.
\end{equation}

The relations (\ref{mam333}) and (\ref{mam334}) imply the equality (\ref{mam331}) for $k=n.$
In view of inequality $\eta<\ye,$ other cases  of the equality (\ref{mam331}) 
can be obtained by using
$$D_{\ye \eta}^{\alpha-k}S_{\ye y}^{n-\alpha}g(\eta)=
(-1)^{n-k}\frac{\partial^{n-k}}{\partial \eta^{n-k}}D_{\ye \eta}^{\alpha-n}S_{\ye y}^{n-\alpha}g(\eta).$$ 
The proof of Lemma 4.2 is complete.

\noindent
{\bf Lemma 4.3.} 
{\it Let $\eta<\ye<y<T,$ $y_{\varepsilon}=y-\varepsilon,$ $n-1<\alpha\leq n,$ $n\in{\mathbb N}.$
Then for the functions $u(\eta)$ and $v(y;\eta)$ such that 
$u(\eta)\in C[0,T],$ $\partial_{0\eta}^{\alpha}u(\eta)\in L[0,y]\cap C(0,\ye),$
$|v(y;\eta)|\leq C(y-\eta)^{\gamma-1},$ $0<\gamma<1,$
the following  estimates 
\begin{equation}\label{lem501}
\left|D_{0\eta}^{\alpha-k}S_{\ye y}^{n-\alpha}v(y;\eta)\right|\leq
C_0(\ye-\eta)^{-\alpha+k-1} \varepsilon^{\gamma}, \quad  k=0,1,...,n,
\end{equation}
\begin{equation}\label{lem50}
\left|\int\limits_{0}^{\ye}\partial_{0\eta}^{\alpha}u(\eta)S_{\ye y}^{n-\alpha}v d\eta
\right|\leq
K_1(\ye,\rho) \varepsilon^{\gamma}+
K_2(\ye,\rho) \varepsilon^{\gamma+n-\alpha-\theta}, \quad  0\leq \theta\leq 1,
\end{equation}
hold, here 
$$K_1(\ye,\rho)=C_1\left[(\ye-\rho)^{-\alpha+k-1}+\ye^{-\alpha+k-1}\right],
\quad
K_2(\ye,\rho)=C_2(\ye-\rho)^{\alpha+\theta-n},$$
$C_0=\frac{C\Gamma(\alpha-k+1)}{\pi\gamma},$
$C_1=\frac{2\Gamma(\alpha-n+1)}{\pi\gamma}CM_0,$
$C_2=\frac{CM_1}{\alpha+\theta-n},$
$M_0=\max\limits_{[0,T]}|u^{(n-1)}(y)|,$
$M_1=\max\limits_{[\rho,\ye]}|\partial_{\rho\eta}^{\alpha}u(\eta)|.$
}

\noindent
{\bf Proof.}
The estimate 
\begin{equation}\label{lem52}
\int\limits_{\ye}^{y}\frac{g(\tau)d\tau}{(\tau-\eta)^{\delta}}\leq
C(\ye-\eta)^{-\delta}\int\limits_{\ye}^{y}(y-\tau)^{\gamma-1}d\tau\leq
\frac{C}{\gamma} 
(\ye-\eta)^{-\delta}(y-\ye)^{\gamma}
\end{equation}
holds.
From (\ref{mam333}) and (\ref{lem52}) we obtain (\ref{lem501}).

Formula integration by parts and relation (\ref{vfh}) transform the following integral
$$\int\limits_{0}^{\ye}\partial_{0\eta}^{\alpha}u(\eta) \cdot S_{\ye y}^{n-\alpha}v(y,\eta) d\eta=
\int\limits_{0}^{\ye}u^{(n)}(\eta) D_{\ye\eta}^{\alpha-n} S_{\ye y}^{n-\alpha}v(y,\eta) d\eta=$$
$$=\left(\int\limits_{0}^{\rho}+\int\limits_{\rho}^{\ye}\right)
u^{(n)}(\eta) D_{\ye\eta}^{\alpha-n} S_{\ye y}^{n-\alpha}v(y,\eta) d\eta=
\left(u^{(n-1)}(\eta)D_{\ye\eta}^{\alpha-n}S_{\ye y}^{n-\alpha}v(y,\eta)\right)\big|_{\eta=0}^{\eta=\rho} + $$
\begin{equation}\label{lem51}
+\int\limits_{0}^{\rho}
u^{(n-1)}(\eta) D_{\ye\eta}^{\alpha-n+1}S_{\ye y}^{n-\alpha}v(y,\eta) d\eta+
\int\limits_{\rho}^{\ye}
\partial_{\rho\eta}^{\alpha}u(\eta) S_{\ye y}^{n-\alpha}v(y,\eta) d\eta.
\end{equation}

From (\ref{lem51}), by using Lemmas 2.1 and 4.2, and the estimate (\ref{lem501}), the 
inequality (\ref{lem50}) follows.
The proof of Lemma 4.3 is complete.

\subsection{The proof of Theorem 3.1}

Integrate equality
$$G{\bf L}u-u{\bf L^*}G=(G\partial_{0\eta}^{\alpha}u-uD_{y\eta}^{\alpha}G) + 
b(G\partial_{0\eta}^{\beta}u-uD_{y\eta}^{\beta}G)-(Gu_{\xi\xi}-uG_{\xi\xi})$$
by the domain 
$\Omega_{\varepsilon,\delta}=\{(\xi,\eta): a_1+\delta_1 <\xi<a_2-\delta_2, 0<\eta<\ye \}.$
By virtue of $Gu_{\xi\xi}-uG_{\xi\xi}=(Gu_{\xi}-uG_{\xi})_{\xi}$ and Lemma 4.1
we obtain
$$\int\limits_0^{y_{\varepsilon}}
\int\limits_{a_1+\delta_1}^{a_2-\delta_2}
[G(x,y;\xi,\eta){\bf L}u(\xi,\eta)-u(\xi,\eta){\bf L^*}G(x,y;\xi,\eta)]d\xi d\eta=$$
$$=\sum\limits_{i=1}^{2}(-1)^{i+1}\int\limits_{0}^{y_{\varepsilon}}
[G(x,y;\bar a_i,\eta)u_{\xi}(\bar a_i,\eta)-
G_{\xi}(x,y;\bar a_i,\eta)u(\bar a_i,\eta)]d\eta+$$
$$+\int\limits_{a_1+\delta_1}^{a_2-\delta_2}
u(\xi,y_{\varepsilon})\left(D_{yy_{\varepsilon}}^{\alpha-1}+bD_{yy_{\varepsilon}}^{\beta-1}\right)
G(x,y;\xi,y_{\varepsilon})d\xi-$$
$$-\int\limits_{a_1+\delta_1}^{a_2-\delta_2}
\left[u\left(D_{\ye \eta}^{\alpha-1}+bD_{\ye \eta}^{\beta-1}\right)G+
(n-1)u_{\eta} D_{\ye\eta}^{\alpha-2}G\right]\Big |_{\eta=0}d\xi+$$
\begin{equation} \label{kmam555}
+\int\limits_{a_1+\delta_1}^{a_2-\delta_2}\!\left[
R_{\alpha}(x,y,y_{\varepsilon},\xi)+bR_{\beta}(x,y,y_{\varepsilon},\xi)
-(n\!-\!1)u_{\eta}(\xi,y_{\varepsilon})D_{yy_{\varepsilon}}^{\alpha-2}G(x,y;\xi,y_{\varepsilon})\right]d\xi,
\end{equation}
where $\bar a_i=a_i-(-1)^i \delta_i,$ $(i=1,2),$
$$R_{\alpha}(x,y,y_{\varepsilon},\xi)= \gamma_n^{\alpha}\int\limits_{0}^{\ye}
\partial_{0\eta}^{\alpha}u(\xi,\eta)S_{\ye y}^{n-\alpha}G(x,y;\xi,\eta)d\eta,$$
$$R_{\beta}(x,y,y_{\varepsilon},\xi)=\gamma_1^{\beta}
\int\limits_{0}^{\ye}\partial_{0\eta}^{\beta}u(\xi,\eta)S_{\ye y}^{1-\beta}G(x,y;\xi,\eta)d\eta.$$

By successively letting  $\delta,$ and $\varepsilon$
tend to zero in relation (\ref{kmam555}), and  by taking into account the relations
(\ref{mam1}), (\ref{oos121}), (\ref{lvh}), (\ref{oos12}),  
$$G(x,y;\xi,\eta)=G(x,y-\eta;\xi,0),$$
$$D_{\ye\eta}^{\mu}\varphi(y-\eta)=D_{\varepsilon z}^{\mu}\varphi(z)\big|_{z=y-\eta},$$
and Lemmas 2.5, 2.7, 2.8 and 4.3
we obtain
(\ref{kkkkk}).
The proof of Theorema 3.1 is complete.

\section{Theorem of existence and uniquiness of solution }
\label{sec:5}

The validity of Theorem 3.2 implies from the Lemmas 2.6 and 2.8.

Let us prove Theorem 3.3.

\noindent
{\bf Proof.}
Theorem 3.1 and Theorem 3.2 imply that every solution of
the problem (\ref{mam1}), (\ref{oos121}), (\ref{zad00}) can be represented in the form (\ref{sp00}). From this representation, we obtain the
solution uniqueness for the investigated problem.
From Lemmas 2.3 and 2.5 we have
$${\bf L}D_{0y}^{\mu}\Gamma(x,y)=0, \quad \mu\in\{\alpha-1,\, \beta-1, \, \alpha-n\}.$$
By virtue of the last relations, Lemmas 2.5, 2.10 and 2.11, 
one can readily find that $u(x, y)$ is a solution of equation (\ref{mam1}) such
that
$\partial_{0y}^{\alpha}u,\, \partial_{0y}^{\beta}u, 
\,\frac{\partial^2}{\partial x^2}u \in C(\Omega).$

Let us prove that conditions (\ref{oos121}) and (\ref{zad00}) is satisfied.
We denote the second and the third integrals on the right side of (\ref{sp00}) by $u_2(x,y)$ and $u_3(x,y)$ respectively,
and by $u_{1,n}(x,y)=u(x,y)-u_2(x,y)-u_3(x,y).$
It follows from Lemma 2.5 and conditions of Theorem 3.3 that
\begin{equation} \label{luf1}
\lim\limits_{y\to 0}\frac{\partial^{k-1}}{\partial y^{k-1}}u_{3}(x,y)=0,  
\quad   k=1,n, \quad a_1\leq x\leq a_2, 
\end{equation}
\begin{equation} \label{luf2}
\lim\limits_{x\to a_i}u_{3}(x,y)=0,  \quad  0\leq y\leq T.
\end{equation}

Taking into account the following estimates
\begin{equation} \label{cl5}
\left|D_{0y}^{\beta-1}\Gamma(x,y)\right|\leq C|x|^{-\theta}y^{\beta\theta}, \quad
\left|D_{0y}^{\alpha-2}\Gamma(x,y)\right|\leq C|x|^{-\theta}y^{\beta\theta-\beta+1}, \quad \theta\geq 0,
\end{equation}
which follow from Lemma 2.5, we obtain
$$\lim\limits_{y\to 0}u_{1,n}(x,y)=\lim\limits_{y\to 0}\int\limits_{a_1}^{a_2}\tau_1(\xi)D_{0y}^{\alpha-1}G(x,y;\xi,0)d\xi.$$
By virtue of the estimate (\ref{kdxdyg}) and the fact that $X_1^m(x,\xi)=0$ only for $m=0$ and $x=\xi$ and
$X_2^m(x,\xi)=0$ 
only in two cases, namely, for 
$m=0$ and $x=\xi=a_1$ and for 
$m=-1$ and $x=\xi=a_2,$
we have the equality
\begin{equation} \label{t32}
\lim\limits_{y\to 0}u_{1,n}(x,y)=\lim\limits_{y\to 0}\int\limits_{a_1}^{a_2}\tau_1(\xi)D_{0y}^{\alpha-1}[\Gamma(x-\xi,y)-\Gamma(x+\xi-2a_1,y)-\Gamma(x+\xi-2a_2,y)]d\xi
\end{equation}
for $a_1\leq x\leq a_2,$ and  equality
$$\lim\limits_{y\to 0}\frac{\partial}{\partial y}u_{1,n}(x,y)=\lim\limits_{y\to 0}\int\limits_{a_1}^{a_2}\tau_1(\xi)\left[D_{0y}^{\alpha}+bD_{0y}^{\beta}\right]\Gamma(x-\xi,y)d\xi+
$$
\begin{equation} \label{t321}
+\lim\limits_{y\to 0}\int\limits_{a_1}^{a_2}\tau_2(\xi)D_{0y}^{\alpha-1}\Gamma(x-\xi,y)d\xi
\end{equation}
for $a_1< x< a_2.$
The first summand in right side of (\ref{t321}) is zero in the condition (\ref{equstk}) \cite{du15-9}.

From the relations (\ref{t32}), (\ref{t321}) and Lemmas 2.5, 2.8 and 2.10, for $x\not=a_i,$ we obtain
\begin{equation} \label{t34}
\lim\limits_{y\to 0}\frac{\partial^{k-1}}{\partial y^{k-1}}u_{1,n}(x,y)=\tau_k(x), \quad k=1,n, \quad a_1<x<a_2.
\end{equation}

In view of the relations 
$X_1^m(x,a_1)=-X_2^m(x,a_1),$  $X_1^{m+1}(x,a_2)=-X_1^m(x,a_2),\,$ 
and
$\frac{\partial X_k^m}{\partial t}=(-1)^k,$  
we obtain
\begin{equation} \label{eq3010}
G_t(x,y;a_i,\eta)=-2\sum\limits_{m=-\infty}^{\infty}
\left(\frac{\partial }{\partial X}\Gamma(X,y-\eta)\right)\big|_{X=X_1^m(x,a_i)}.
\end{equation}
Since $X_1^m(x,a_i)=0$ only for $m=0$ and $x=a_i$ we have
\begin{equation} \label{a}
\lim\limits_{y\to 0}\int\limits_{0}^{y}\varphi_i(\eta)G_{\xi}(x,y;a_i,\eta)d\eta=
-2\lim\limits_{y\to 0}\int\limits_{0}^{y}\varphi_i(\eta)
\frac{\partial }{\partial x}\Gamma(x-a_i,y-\eta)d\eta, 
\end{equation}
for $ a_1\leq x\leq a_2.$
 
It follows from last relation and the estimate 
$$\left|\frac{\partial^{k-1}}{\partial y^{k-1}}\frac{\partial }{\partial x}\Gamma(x,y)\right|\leq C|x|^{-\theta}y^{\beta\theta-k}, \quad \theta\geq 0,$$
which holds by virtue of Lemma 2.5, that the relation
\begin{equation} \label{dyu2}
\lim\limits_{y\to 0}\frac{\partial^{k-1} }{\partial y^{k-1}}u_2(x,y)=0, \quad a_1<x<a_2
\end{equation}
holds for each $x\not=a_i.$

By taking into account  the relations
$X_1^m(a_1,\xi)=-X_2^{-m}(a_1,\xi),$ $X_1^m(a_2,\xi)=-X_2^{-m-1}(a_2,\xi),$ 
and  that
$X_1^m(a_1,\xi)\not=0,$ for $m\not=0,$   and
$X_1^m(a_2,\xi)\not=0,$ for $m\not=-1,$  we obtain
$$\lim\limits_{x\to a_i}u_{1,n}(x,y)=\lim\limits_{x\to a_i}
\sum\limits_{k=1}^{n}\int\limits_{a_1}^{a_2}\tau_k(\xi)\left[D_{0y}^{\alpha-k}+(2-k)bD_{0y}^{\beta-k}\right]
\Phi(x,\xi,y)d\xi,
$$
where
$\Phi(x,\xi,y)=\Gamma(x-\xi,y)-\Gamma(x+\xi-2a_i,y).$
Since, by virtue of (\ref{cl5}), the following inclusion 
$$D_{0y}^{\alpha-2}\Phi(x,\xi,y), D_{0y}^{\beta-1}\Phi(x,\xi,y)\in C([a_1,a_2]\times\overline{\Omega}),$$
holds, then
\begin{equation} \label{b}
\lim\limits_{x\to a_i}u_{1,n}(x,y)=\lim\limits_{x\to a_i}
\int\limits_{a_1}^{a_2}\tau_1(\xi)D_{0y}^{\alpha-1}\Phi(x,\xi,y)d\xi, \quad a_1\leq x\leq a_2. 
\end{equation}
And since 
$$D_{0y}^{\alpha-1}\Phi(x,\xi,y)\in C((a_1,a_2)\times\overline{\Omega}),$$
then
\begin{equation} \label{c}
\lim\limits_{x\to a_i}u_{1,n}(x,y)=0, \quad a_1< x< a_2. 
\end{equation}

From the relation (\ref{eq3010}), Lemmas 2.5, 2.7 and 2.9 we have
\begin{equation} \label{d}
\lim\limits_{x\to a_i}\int\limits_0^yG_{\xi}(x,y;a_i,\eta)\varphi_i(\eta)d\eta=(-1)^i\varphi_i(y), \quad (i=1,2), \quad y\not=0,
\end{equation}
\begin{equation} \label{e}
\lim\limits_{x\to a_i}\int\limits_0^yG_{\xi}(x,y;a_{3-i},\eta)\varphi_{3-i}(\eta)d\eta=0, \quad (i=1,2), \quad y\not=0.
\end{equation}

Formulas (\ref{luf1}), (\ref{luf2}), (\ref{t34}), (\ref{dyu2}), (\ref{c}), (\ref{d}) and (\ref{e}) imply
that $u(x,y)$ is a solution of problem (\ref{mam1}), (\ref{oos121}), (\ref{zad00}) such that
$u(x,y)\in C(\overline{\Omega}\setminus \{(a_i,0), i=1,2\}).$

Let us prove the continuity of the function $u(x,y)$ at the points $(a_i,0), i=1,2.$
It follows from the relations (\ref{luf1}), (\ref{luf2}), (\ref{t32}), (\ref{a}) and 
$$2\frac{\partial }{\partial x}\Gamma(x,y)=
\int\limits_{|x|}^{\infty}g(y,\tau)\frac{\partial}{\partial x}h_0(x,\tau)d\tau
-{\rm sgn}(x)g(y,|x|),$$
and the estimate \cite{du16-6}
$$\left|\int_{|x|}^{\infty}g(y,\tau)
\frac{\partial}{\partial x}h_0(x,\tau)d\tau\right|\leq 
C |x|^{1-\theta}y^{\beta+\beta\theta-1}, \quad \theta\geq 0, \quad |x|\geq 0,
$$
 that
$$\lim\limits_{x\to a_i \atop y\to 0}u(x,y)=\lim\limits_{x\to a_i \atop y\to 0}\int\limits_{0}^{y}g(y-\eta,|x-a_i|)\varphi_i(\eta)d\eta+
$$
$$ +\lim\limits_{x\to a_i \atop y\to 0}
\int\limits_{a_1}^{a_2}\tau_1(\xi)D_{0y}^{\alpha-1}[g(y,|x-\xi|)-g(y,|x+\xi-2a_i|)]d\xi.$$
Similarly, like in \cite{du16-6} we can prove that  
$$\lim\limits_{x\to a_i \atop y\to 0}u(x,y)= \tau_1(a_i)+ \phi(-\beta,1;-c)[\varphi_i(0)-\tau_1(a_i)],$$
where $c=\lim\limits_{x\to a_i \atop y\to 0}|x-a_i|y^{-\beta}.$
Hence the limit $\lim\limits_{x\to a_i \atop y\to 0}u(x,y)$ is independent of $c$ under the condition (\ref{usog}).
The proof of Theorem 3.3 is complete.





 \bigskip \smallskip

 \it

\noindent
Institute of Applied Mathematics and Avtomation \\
"Shortanov" $\,$ Str., 89 A, \\
Nal'chik -- 360000, RUSSIA  \\[4pt]
e-mail: mamchuev@rambler.ru


\begin{thebibliography}{99}
\normalsize


\bibitem{n} 
A. M. Nakhushev,  
\emph{Fractional Calculus and Its Application}.
Fizmatlit, Moscow, (2003) (In Russian).

\bibitem{Hilfer} 
R. Hilfer,  
\emph{Applications of fractional calculus in physics}. 
World Scientific, New Jersey,  (2000). 

\bibitem{Uchaikin} 
V.V. Uchaikin,
\emph{The method of fractional derivatives}.
Artishok, Ulianovsk, 2008. 512 p. (In Russian). 


\bibitem{Cascaval-2002}
R. C. Cascaval, E. C. Eckstein, C. L. Frota, and J. A. Goldstein,
Fractional telegraph equations.
\emph{J. Math. Anal. Appl.} \textbf{276}, No 1 (2002), 145--159.


\bibitem{Orsinger-2003}
E. Orsinger, and X. Zhao, 
The space-fractional telegraph equation and the related fractional telegraph process. 
\emph{Chinese Ann. Math. Ser. B.} \textbf{24}, No 1 (2003), 45--56.

\bibitem{Orsinger-2004}
E. Orsinger, and L. Beghin,
Time-fractional telegraph equations and telegraph processes with brownian time. 
\emph{Probab. Theory Related Fields.} \textbf{128}, No. 1 (2004), 141--160.

\bibitem{Rekh2016}
S. Sh. Rekhviashvili, Murat O. Mamchuev, and Mukhtar O. Mamchuev, 
Model of diffusion-drift charge carrier transport in layers with a fractal structure.
\emph{Physics of the Solid State.} \textbf{58}, No 4 (2016), 789--792.



\bibitem{Liu-2003} F. Liu, V. V. Anh, I. Turner and P. Zhuang, 
Time fractional advection-dispersion equation. 
\emph{J. Appl. Math. Computing.} \textbf{13}, (2003), 223--245.

\bibitem{Huang-2005} F. Huang and F. Liu, 
The time fractional diffusion equation and the advection-dispersion equation. 
\emph{The A N Z I A M Journal.} \textbf{46}, (2005), 317--330.


\bibitem{izvkbnc04} 
M. O. Mamchuev,  
General representation of a solution of a fractional diffusion equation 
with constant coefficients in a rectangular domain. 
\emph{Izv. Kabardino-Balkarsk. Nauch. Ts. RAN.} \textbf{12}, No 2 (2004), 116--118 (In Russian).

\bibitem{daman05} 
M. O. Mamchuev,  
Boundary value problems for a fractional diffusion equation with constant coefficients. 
\emph{Dokl. Adyg. (Cherkessk.) Mezhdunar. Akad. Nauk.} \textbf{7}, No 2 (2005), 38--45 (In Russian).

\bibitem{du15-5} 
M. O. Mamchuev,  
Fundamental solution of a loaded second-order parabolic equation with constant coefficients. 
\emph{Differential Equations.} \textbf{51}, No 5 (2015), 620--629.

\bibitem{du15-9} 
M. O. Mamchuev,  
Modified Cauchy problem for a loaded second-order parabolic equation with constant coefficients. 
\emph{Differential Equations.} \textbf{51}, No 9 (2015), 1137--1144.


\bibitem{Mamchu2016arx}
Murat O. Mamchuev,
Solutions of the main boundary value problems for the time-fractional telegraph equation by the Green function method.
arXiv:1608.06852.  22 pages.




\bibitem{Ps-2003} 
A. V. Pskhu, Solution of boundary value problems for the fractional diffusion
equation by the Green function method. 
\emph{Differential Equations.} \textbf{39}, No 10 (2003), 1509--1513.

\bibitem{PsMon}
A. V. Pskhu,  
\emph{Fractional Partial Differential Equations}. Nauka, Moscow (2005) (In Russian). 





\bibitem{Caputo}
{\it Caputo M.} Elasticita e Dissipazione. Zanichelli, Bologna, 1969. 
(in Italian).


\bibitem{Kilbas-2006}
A.A. Kilbas, N.H.M. Srivastava and J.J. Trujillo, 
\emph{Theory and applications of fractional differential equation}. 
Elsevier, Amsterdam, 2006.

\bibitem{Wright-1933}
E. M. Wright, On the coeffcients of power series having exponential singularities. 
\emph{J. London Math. Soc.} \textbf{8}, (1933), 71--79.


\bibitem{du16-6} 
M. O. Mamchuev,  
Solutions of the main boundary value problems for a loaded second-order parabolic equation with constant coefficients. 
\emph{Differential Equations.} \textbf{52}, No 6 (2016), 789--797.










\end{thebibliography}
\end{document}